\documentclass[12pt,a4paper]{article}
\usepackage{amsfonts}

\usepackage{amssymb}
\usepackage{amscd}
\usepackage{amsmath}

\begin{document}

\begin{center}
{\huge Algebraic Smooth Structures}{\large \ 1}{\Large \\[0pt]
}\vspace{0.4cm}

{\small {Shafiei Deh Abad A.\\[0pt]
School of mathematics, Statistics and Computer Science, University of
Tehran, Tehran, Iran}\\[0pt]
\vspace{0.1cm}E-mail : shafiei@khayam.ut.ac.ir}

{\small \thinspace\ }

{\small {\large Abstract}\vspace{0.1cm} }
\end{center}

{\small In this paper which is the first of a series of papers on smooth
structures, the concepts of }$C${\small -structures and smooth structures
are introduced and studied. The notion of smooth structure on semi-integral
domains is given. It is shown that each semi-integral domain which is not a
field, admits a unique smooth structure and a large class of non-polynomial
smooth functions on some semi-integral domains is constructed. A smooth
function from }$\
\mathbb{Z}
${\small -\{0\} into }$%
\mathbb{Z}
${\small \ is given which does not extend to a smooth function on}$\
\mathbb{Z}
\ ${\small . No concept from topology is used. As an application, it is
shown that:}

{\small Theorem - Let }$M${\small \ and }$N${\small \ be finite dimensional
smooth manifolds. Assume that }$\varphi :C^{\infty }(N)\rightarrow C^{\infty
}(M)${\small \ is a homomorphism of }$\
\mathbb{R}
\ ${\small -algebras. Then, there exists exactly one smooth mapping }$\phi
:M\rightarrow N${\small \ such that }$\varphi =\phi ^{\ast }${\small . \ \ }$%
\ \ \ \ \ $

\bigskip

\section{\protect\Large \ \ Introduction}

\bigskip

Differential calculus is a powerful technic in mathematics. To name a few,
partial differential equations, differential and analytic geometry are all
based on the notion of differentiation, and most of the applications of
mathematics in other sciences employ this concept. The whole theory of
differentiation is based on topology and the fact that the underlying rings
are fields. A construction of a satisfactory theory of differentiation on $%
\mathbb{R}
$ or $%
\mathbb{C}
$ without any use of topology is invisible. On the other hand, one can
redefine the notion of derivative using only the topology and the ring
structure of the fields. More precisely, it is not difficult to see that the
following definition of the derivative of a function on $%
\mathbb{R}
$ is equivalent to the usual one.

{\bf Definition 1.1 \ }Let $f$ be a real valued function defined on a
neighborhood of $\lambda \in
\mathbb{R}
$. Then $f$ is differentiable at $\lambda $ and its derivative at $\lambda $
is $\alpha $, if and only if there exists a real valued function g defined
on V, a neighborhood of $\lambda $, and continuous at this point, such that

\begin{center}
\[
g(\lambda )=\alpha\quad and\quad f(t)=f(\lambda )+(t-\lambda )g(t)
\quad for\; all\quad t\in V.
\]
\end{center}

It is clear that the function g with the above properties is unique and the
only use of topology here is to characterize g($\lambda $) uniquely. Clearly
if we can avoid this usage we can define derivative of functions without
ambiguity. This was what we have done some years ago in [3] where we define
the derivative of functions on integral domains which are not fields,
without using any concept from topology. Such a derivative has properties
similar to the properties of usual derivative of real functions and gives
rise to very interesting problems. In this paper the definition of smooth
structure is extended to the functions on semi-integral domains. The smooth
structure in this general case have very nice properties$.$ For example, let
$R$ be a proper semi-integral domain. Then:

1) The $R$-algebra of smooth functions on $R$ is itself a proper
semi-integral domain.

2) The map $\varphi :C^{\infty }(R^{2})\longrightarrow C^{\infty
}(R,C^{\infty }(R))$ given by $\varphi (f)(x)(y)=f(x,y)$ is an isomorphism.

3) The $R$-algebra of smooth functions on $R^{n}$, $n\succeq 1$
is uniquely determined by $C^{\infty }(R)$ without any use of
topology.

On the other hand, differential calculus on rings has very
peculiar properties. For example, on some integral domains $R$
there exist

i) non-constant smooth functions with derivatives identically zero.

ii) functions $f$ and $0\neq \lambda \in R$ such that $\lambda f$
is smooth but $f$ is not.

Our next goal is to define smooth structures on modules. The present paper
is devoted to "smooth structures". It serves as a foundation for the
subject. In addition we will provide interesting applications.

The work will be continued in the forthcoming papers. We will prove that any
projective modules over a semi-integral domain which is not a field\ admits
a unique smooth structure, and differential calculus on finitely generated
projective modules almost always has the same properties as the usual
differential calculus on $%
\mathbb{R}
^{{\it n}}.$ The proof of Propositions 2 , 3 and 4 above will be
given in appropriate places in forthcoming papers.$\bigskip $\
\section{$\protect\bigskip ${\protect\Large \ \ Conventions }}

\bigskip In what follows $R$ denotes a commutative ring with identity
element 1. By an $R-$algebra we mean an associative and commutative $R-$%
algebra with an identity. Furthermore, for every $R-$algebra with an
identity element $e$, we assume that $\left\{ \lambda \in R\mid \lambda
e=0\right\} =\left\{ 0\right\} .$ All algebra homomorphisms are assumed to
preserve the identity element. Any sub-algebra contains the identity
element. Finally, by a sub-module of an $R-$algebra ${\cal A}$ we mean a
sub-module of the $R-$module ${\cal A}$.

Let ${\cal {A}}$ be an $R$-algebra with the identity e. Assume that $A$ and $%
B$ are subsets of ${\cal {A}}$. The sub-module of ${\cal {A}}$ generated by
the set $\{{xy\mid x\in \ A,y\in \ B}\}$ will be denoted by $A.B$. For each $%
n\in N^{\ast }$, let $^{n}A=\{{x_{1}x_{2}\ldots x_{n}\mid x_{i}\in A}\}$.
The sub-module of ${\cal {A}}$ (resp. the ideal of ${\cal {A}}$) generated
by $^{n}A$ will be denoted by $A^{n}$ (resp. by $\underline{A}^{n})$. We use
the convention ${\cal {A}}^{o}$ = Re. Moreover if $A$ is a sub-module of $%
{\cal {A}},A^{1}=A$

\section{\protect\bigskip\ {\protect\large {\bf \ Algebras of type C}}}

\bigskip Let ${\cal {A}}$ be an $R$-algebra with the identity e. An ideal $%
\Lambda $ of ${\cal {A}}$ is called {\it an ideal of type} $C$ if ${\cal {A}}
$ = Re $\oplus \Lambda $. We identify $Re$ with $R$ and denote the
projection on $R=Re$ (resp. on $\Lambda )$ by $\pi _{\Lambda }^{0}$ (resp.
by $\pi _{\Lambda })$. The set of all ideals of type $C$ in ${\cal {A}}$
will be denoted by ${\cal {A}}^{c}$. We say that ${\cal {A}}$ is an $R$-{\it %
algebra of type} $C$ if ${\cal {A}}^{c}\neq \phi $. Let $\Lambda \in {\cal {A%
}}^{{c}}$. Then $\pi _{\Lambda }^{0}$ is clearly a character of ${\cal {A}}$%
. On the other hand, if $\alpha $ is a character of ${\cal {A}}$, then ker$%
\alpha \in {\cal {A}}^{{c}}$. The set of all characters of ${\cal {A}}$ will
be denoted by ${\cal {A}}^{\gamma }$. The mapping ${\cal {A}}^{\gamma
}\rightarrow {\cal {A}}^{{c}}$ given by $\alpha \mapsto $ ker$\alpha \ $is$\
$clearly bijective.

\ {\bf Lemma 3.1. \ }Let $R$ be an integral domain and let ${\cal A}$ be an $%
R-$algebra of type $C$. Then each $\Lambda \in {\cal A}^{c}$ is a prime
ideal.

\ {\bf Proof. \ }Let $x,y$ be in ${\cal A}$ and let $xy\in \Lambda \in {\cal %
A}^{c}$. Then $\pi _{\Lambda }^{0}(xy)=$ $\pi _{\Lambda }^{0}(x)\pi
_{\Lambda }^{0}(y)=0.$ Since $R$ is an integral domain $\pi _{\Lambda
}^{0}(x)=0$ or $\pi _{\Lambda }^{0}(y)=0.$ Therefore, $x\in \Lambda $ or $%
y\in \Lambda .$

The proof of the following lemma is straightforward.

\hspace{-0.7cm} \ \ \ {\bf Lemma 3.2. }Let ${\cal {A}}$ and ${\cal {A}}%
^{\prime }$ be $R$-algebras and let $\varphi :{\cal {A}}\rightarrow {\cal {A}%
}^{\prime }$ be an $R$-algebra homomorphism. Assume that ${\cal {A}}^{\prime
}$ is an $R$-algebra of type $C$. Then ${\cal {A}}$ is also an $R$-algebra
of type $C$. Moreover, if $\Lambda ^{\prime }\in {\cal {A}}^{\prime c}$ and $%
\alpha ^{\prime }\in {\cal {A}}^{\prime \gamma }$, then $\varphi
^{-1}(\Lambda )\in {\cal {A}}^{{c}}$ and $\varphi ^{\ast }(\alpha ^{\prime
})\in {\cal {A}}^{\gamma }$.

Let ${\cal {A}},{\cal {A}}^{\prime }$ and $\varphi $ be as above. The
mapping ${\cal {A}}^{\prime c}\rightarrow {\cal {A}}^{c}$ given by $\Lambda
^{\prime }\mapsto \varphi ^{-1}(\Lambda ^{\prime })$ will be denoted by $%
\underline{\varphi }^{{\tiny {c}}}$. We say that $\underline{\varphi }^{%
{\tiny {c}}}$ is {\it induced by} $\varphi $.

\hspace{-0.7cm} \ \ \ \ {\bf Lemma\ 3.3. }\ Let ${\cal {A}}$ and ${\cal {A}}%
^{\prime }$ be two $R$-algebras of type $C$. Assume that $\varphi :{\cal {A}}%
\rightarrow {\cal {A}}^{\prime }$ is an $R$-algebra homomorphism. Let $%
\Lambda ^{\prime }\in {\cal {A}}^{\prime {\tiny {c}}}$ and
$\Lambda =\varphi ^{-1}({\Lambda }^{\prime })$. Then, $\pi
_{\Lambda }^{0}=\pi _{\Lambda^{\prime} }^{0}\circ \varphi $ and
$\pi _{\Lambda ^{\prime }}\circ \varphi =\varphi \circ \pi
_{\Lambda }$.

{\bf Proof. }\ Let $y\in {\cal {A}}$. Then $y=\pi _{\Lambda }^{0}(y)e+\pi
_{\Lambda }(y)$ and $\varphi (y)=\pi _{\Lambda ^{\prime }}^{0}(\varphi
(y))e^{\prime }+\pi _{\Lambda ^{\prime }}(\varphi (y))$. Hence, $\pi
_{\Lambda ^{\prime }}^{0}(\varphi (y))e^{\prime }+\pi _{\Lambda ^{\prime
}}(\varphi (y))=\varphi (y)=\varphi \lbrack \pi _{\Lambda }^{0}(y)e+\pi
_{\Lambda }(y)]=\pi _{\Lambda }^{0}(y)e^{\prime }+\varphi (\pi _{\Lambda
}(y))$. Or, $(\pi _{\Lambda ^{\prime }}^{0}\circ \varphi (y)-\pi _{\Lambda
}^{0}(y))e^{\prime }+(\pi _{\Lambda ^{\prime }}\circ \varphi (y)-\varphi
\circ \pi _{\Lambda }(y))=0$. But $\varphi (\Lambda )\subset \Lambda
^{\prime }$ and $\Lambda ^{\prime }\in {\cal {A}}^{\prime {\tiny {c}}}$.
Therefore, $\pi _{\Lambda ^{\prime }}^{0}\circ \varphi (y)=\pi _{\Lambda
}^{0}(y)$ and $\pi _{\Lambda ^{\prime }}\circ \varphi (y)=\varphi \circ \pi
_{\Lambda }(y)$. Since $y\in {\cal {A}}$ is arbitrary, $\pi _{\Lambda
^{\prime }}^{0}\circ \varphi =\pi _{\Lambda }^{0}$ and $\pi _{\Lambda
^{\prime }}\circ \varphi =\varphi \circ \pi _{\Lambda }$.$\blacksquare $

A sub-module $A$ of ${\cal {A}}$ is called a {\it sub-module of type} $C$ if

1) $\underline{A}\in {\cal {A}}^{c}$.

2) If $B$ is a sub-module of $A$ and $\underline{B}$\ $\in {\cal {A}}^{{c}},$
then $A=B.$

The set of all sub-modules of type $C$ of ${\cal {A}}$ will be denoted by $C(%
{\cal {A}})$. A sub-module $A\in C({\cal {A}})$ is called a {\it
sub-module of type} $D$ if for each $k\geq 1,$ $A^{k}\cap
\underline{A}^{k+1}=\{{0}\}$. The set of all sub-modules of type
$D$ will be denoted by $D({\cal {A}})$. If $A\in C({\cal {A}})$
is of type $D$, then $\underline{A}$ is called an ideal of type
$D$. The set of all ideals of type $D$ will be denoted by ${\cal
{A}}^{d}$.
Let $A\in C({\cal {A}})$ and $k\in $ $%
\mathbb{N}
$. Then, $\underline{A}^{k}=A^{k}.{\cal {A}}=A^{k}.(A^{0}+\underline{A}%
)=A^{k}.(A^{0}+A.{\cal {A}})=A^{k}+A^{k+1}.{\cal {A}}=A^{k}+\underline{A}%
^{k+1}$. Therefore, we have the following simple lemma.

\hspace{-0.7cm} \ \ \ \ \ {\bf Lemma 3.4. \ }Let ${\cal {A}}$ be an $R$%
-algebra and $n\in
\mathbb{N}
.$ Then:

i) if $A \in C({\cal {A}})$, then ${\cal {A}} = \sum_{k=0} ^{n} A^{k} +
\underline A^{n+1}$.

ii) if $A\in D({\cal {A}})$ $,$then $\ {\cal {A}}=\oplus
_{k=0}^{n}A^{k}\oplus \underline{A}^{n+1}.$

Let $A$ be a sub-module of type $D$. By the above lemma for $n\geq 1$ , $%
{\cal {A}}=\oplus _{k=0}^{n}A^{k}\oplus \underline{A}^{n+1}$. The projection
of ${\cal {A}}$ onto $A^{k}$ (resp. onto $\underline{A}^{n+1})$ in the above
decomposition will be denoted by $\alpha _{k}$ (resp. by $\underline{\alpha }%
_{n+1})$ and $\pi _{\underline{A}}^{0}:{\cal A\rightarrow }R\cdot
e\backsimeq R$ will be denoted by $\alpha _{0}.$

{\bf Lemma 3.5. \ }Let{\bf \ }$\varphi :{\cal A}\rightarrow {\cal
A}^{\prime }$ be injective. Assume that $A\in C\left( {\cal
A}\right) $ and $\varphi \left( A\right)\subset{A^{\prime}} \in
D\left( {\cal A}^{\prime }\right) .$ Then, $A\in D\left( {\cal
A}\right) $.

{\bf Proof. \ }Let $A\in C\left( {\cal A}\right) $ and $\varphi \left(
A\right) \in D\left( {\cal A}^{\prime }\right) $. Then, for each $0$ $\leq
n, $ $\left( \varphi \left( A\right) \right) ^{n}\cap \left( \underline{%
\varphi \left( A\right) }\right) ^{n+1}=\{0\}.$ Since $\varphi $ is
injective, $A^{n}\cap \underline{A}^{n+1}=\left( \varphi ^{-1}\left( \varphi
\left( A\right) \right) \right) ^{n}\cap \left( \varphi ^{-1}(\varphi
\underline{\left( A\right) }\right) ^{n+1}\subset \varphi ^{-1}(\left(
\varphi \left( A\right) \right) ^{n}\cap \left( \underline{\varphi \left(
A\right) }\right) ^{n+1})=\{0\}.$ Therefore, $A\in D\left( {\cal A}\right) $.$%
\blacksquare $

\bigskip Let \ $\Lambda $ and $\Lambda ^{\prime }$ be in ${\cal A}^{c}.$
Then, for each $y\in {\cal A}$ we have $\pi _{\Lambda ^{\prime }}\left(
y\right) =\pi _{\Lambda ^{\prime }}(\pi _{\Lambda }^{0}\left( y\right) +\pi
_{\Lambda }\left( y\right) )=\pi _{\Lambda ^{\prime }}\circ \pi _{\Lambda
}\left( y\right) .$ Thus, $\pi _{\Lambda ^{\prime }}\circ \pi _{\Lambda
}=\pi _{\Lambda ^{\prime }},$ and for every $x\in \Lambda ^{\prime }$ we
have $x=\pi _{\Lambda ^{\prime }}\left( x\right) =\pi _{\Lambda ^{\prime
}}\circ \pi _{\Lambda }\left( x\right) .$ Therefore, $\pi _{\Lambda ^{\prime
}}:\Lambda \rightarrow \Lambda ^{\prime }$ is an isomorphism of $R-$modules
with inverse $\pi _{\Lambda }.$

Let $A$ and $B$ be in $C({\cal {A}})$. We say that $A$ is{\it \ equivalent to%
} $B$ and write $A\thickapprox B,$ if
$\underline{A}=\underline{B}$. The sub-modules $A$ and $B$ are
called{\it \ strongly compatible with} {\it each other} if:

\[
\ \ A\thickapprox \pi _{\underline{A}}(B)\quad { and }\quad B\thickapprox \pi _{%
\underline{B}}(A).
\]
This will be written as $A\sim B.$ Two sub-modules $A$ and $B$ are called
{\it compatible} if there exists a finite sequence $A=C_{0}\sim C_{1}\sim
C_{2}\sim ...\sim C_{k}=B.$

{\bf Lemma\ 3.6.}\ \ Let $A,B\in D({\cal {A}})$ be two equivalent
sub-modules. Then

i) $\ A$ is strongly compatible with $B$,

\bigskip ii) \ $\alpha _{1}:B\rightarrow A$ is an isomorphism with inverse $%
\beta _{1}$ .\

{\bf Proof.}\ i) This is clear.

ii) \ By (3.4.ii) $Re\oplus A\oplus \underline{A}^{2}={\cal
{A}}=Re\oplus B\oplus \underline{B}^{2}$. Assume that $y\in A$.
The equality $y=\beta _{1}(y)+\underline{\beta }_{2}(y)$ implies
that $y=\alpha _{1}(y)=\alpha _{1}(\beta _{1}(y))$. Thus $\alpha
_{1}\circ \beta _{1}=id\mid _{A}.$ In the same way one sees that
\ $\beta _{1}\circ \alpha _{1}=id\mid _{B}.$ Therefore, $\alpha
_{1}$ is an isomorphism with inverse $\beta _{1}.$

\bigskip Now it is clear that we have the following proposition.

{\bf Proposition 3.7 \ }The relation $A${\it \ is compatible with }$B$ is an
equivalence relation.

Let $M$ be a set. Assume that ${\cal {A}}$ is a subalgebra of $R^{M}$. Let $%
\omega \in M$. Clearly the set of all elements of ${\cal {A}}$ which are
zero at $\omega $ is an ideal of type $C$. It will be denoted by $I_{\omega
} $. If the elements of ${\cal {A}}$ separate the points of $M$, then the
mapping $\sigma :M\rightarrow {\cal {A}}^{c}$ given by $\omega \rightarrow
I_{\omega }$ is injective. In this situation we identify $M$ with its image
under the above $\ $mapping.

\section{{\protect\large {\bf \ \ }}{\protect\Large C-Structures}}

\bigskip

\ \ \ \ A$\ C$-{\it pair} (resp. A$\ D$-{\it pair})  over $R$ is a pair $({\cal {A}},\Sigma )$, where $%
{\cal {A}}$ is an $R$-algebra of type $C,$ and $\Sigma $ is a
non-empty subset of $A$\ C$-{\it pair}$ (resp. $A$\ D$-{\it
pair})$. Let $\Lambda \in \Sigma .$ The set of all $A\in D\left(
{\cal A}\right) $ such that $\underline{A}=\Lambda ,$ will be\
denoted by $\Lambda ^{\circ },$ and $\cup _{\Lambda \in \Sigma
}\Lambda ^{\circ }$ will be denoted by $\overline{\Sigma }.$

Let $({\cal A}$, $\Sigma )$ and $({\cal A}^{\prime }$, $\Sigma ^{\prime })$
be two ${\it C-}$pairs over $R$. A $C-$homomorphism $\varphi :({\cal A}$, $%
\Sigma )\rightarrow ({\cal A}^{\prime }$, $\Sigma ^{\prime })$ is a
homomorphism of $R-$algebras $\varphi :{\cal A\rightarrow A}^{\prime }$ such
that for all $\Lambda ^{\prime }\in \Sigma ^{\prime },$ $\varphi
^{-1}(\Lambda ^{\prime })\in \Sigma .$ Clearly, the restriction of $%
\underline{\varphi ^{c}}$ to $\Sigma ^{\prime }$ is a map from $\Sigma
^{\prime }$ into $\Sigma .$ This map will be denoted by \underline{$\varphi $%
}$.$ It is also clear that the map {\it id}$_{{\cal A}}:({\cal A}$, $\Sigma
)\rightarrow ({\cal A}$, $\Sigma )$ is a $C-$homomorphism. Let $({\cal A}%
_{i},\Sigma _{i})$,$\ i=1,2,3$ be $C-$pairs over $R$. Assume that $\varphi
_{1}:({\cal A}_{1}$, $\Sigma _{1})\rightarrow ({\cal A}_{2}$, $\Sigma _{2})$
and $\varphi _{2}:({\cal A}_{2},\Sigma _{2})\rightarrow ({\cal A}_{3},\Sigma
_{3})$ are $C-$homomorphisms. Then, clearly $\varphi _{2}\circ \varphi _{1}$
is a $C-$homomorphism, and $\underline{\varphi _{2}\circ \varphi _{1}}=%
\underline{\varphi _{1}}\circ \underline{\varphi _{2}}.$

The above observations can be summarized in the following.

{\bf Proposition\ 4.1. \ }The class of all $C-$pairs over $R$ together with $%
C-$homomorphisms between them form a category. This category will be denoted
by $R-CP.$

A $C-$homomorphism $\varphi:({\cal A},\Sigma){\cal \rightarrow (A}%
^{\prime},\Sigma^{\prime})$ is called injective (resp.
surjective, resp. bijective) if $\varphi :{\cal A\rightarrow
A}^{\prime }$ is injective (resp. surjective, resp. bijective).

A $C-$pair $({\cal A},\Sigma )$ over $R$ is called

i){\it \ separated} if $\cap _{\Lambda \in \Sigma }\Lambda =\{0\},$

ii) {\it analytic}\ if for each $\Lambda \in \Sigma $, $\cap
_{n=1}^{\infty }\Lambda ^{n}=0,$

iii) {\it of polynomial type }if for each $\Lambda \in \Sigma $ there exists
$A\in D\left( {\cal A}\right) $ such that $\underline{A}=\Lambda $, and each
element $y\in {\cal A}$ can be written as $y=\Sigma _{n=0}^{k}y_{n}$, where $%
y_{n}\in A^{n}.$

Assume that the $C-$pair $({\cal A},\Sigma )$ is separated. Then the
mapping\ $\ \theta :{\cal A}\rightarrow R^{\Sigma }$ given by $y\mapsto
(\Lambda :\mapsto \pi _{\Lambda }^{0}(y))$ is clearly an injective
homomorphism of $R-$algebras. In this situation we identify ${\cal A}$ with
its image under the mapping $\theta .$

{\bf Lemma\ 4.2. }\ Let $\varphi :({\cal A},\Sigma )\rightarrow ({\cal A}%
^{\prime },\Sigma ^{\prime })$ be a $C-$homomorphism. Then

i) If $\varphi$ is injective and $({\cal A}^{\prime},\Sigma^{\prime})$ is
separated, then $({\cal A}$,$\Sigma)$ is separated.

ii) If $({\cal A},\Sigma)$ is separated and \underline{$\varphi$} is
surjective, then $\varphi$ is injective.

iii) If $\varphi $ is surjective then \underline{$\varphi $} is injective.

\bigskip {\bf Proof. }i) We have{\bf \ }$\cap _{\Lambda \in \Sigma }\Lambda
\subset \cap _{\Lambda ^{\prime }\in \Sigma ^{\prime }}\varphi ^{-1}(\Lambda
^{\prime })=\varphi ^{-1}(\cap _{\Lambda ^{\prime }\in \Sigma ^{\prime
}}\Lambda ^{\prime })=\varphi ^{-1}(\{0\})=\{0\}.$ Therefore, $({\cal A}%
,\Sigma )$ is separated.

ii) Let $y\in {\cal A}.$ Assume that $\varphi \left( y\right) =0.$ By Lemma
3.2 for each $\Lambda ^{\prime }\in \Sigma ^{\prime },$ \ $\pi _{\underline{%
\varphi }(\Lambda ^{\prime })}^{0}(y)=\pi _{\Lambda ^{\prime }}^{0}(\varphi
(y))=0.$ Since $\underline{\varphi }$:$\Sigma ^{\prime }\rightarrow \Sigma $
is surjective, for all $\Lambda \in \Sigma ,$ $\pi _{\Lambda }^{0}(y)=0.$
But $({\cal A},\Sigma )$ is separated. Therefore, $y=0$ and $\varphi $ is
injective.

iii) Let $\Lambda ^{\prime },\Lambda ^{\prime \prime }\in \Sigma
^{\prime }$ and $\varphi ^{-1}(\Lambda ^{\prime })=\Lambda
=\varphi ^{-1}(\Lambda ^{\prime \prime })$. Since $\varphi $ is
surjective $\Lambda ^{\prime \prime }=\varphi \left( \varphi
^{-1}\left( \Lambda ^{\prime \prime }\right) \right) =\varphi
\left( \Lambda \right) =\varphi \left( \varphi ^{-1}\left(
\Lambda ^{\prime }\right) \right) =\Lambda ^{\prime }.$
Therefore, $\underline{\varphi} $ is injective. $\blacksquare $\ \
\ \ \ \ \ \ \ $\ \ \ \ \ \ \ \ \ \ \ \ \ \ \ \ \ \ \ \ \ \ \ \ \ \
\ \ \ \ \ \ \ \ \ \ \ \ \ \ \ \ \ \ \ \ $\ \ \ \ \ \ \ \ \ \ \ \
\ \ \ \ \ \ \ \ \ \ \ \ \ \ \ \ \ \ \ \ \ \ \ \ \ \ \ \ \ \ \ \ \
\ \ \ \ \ \ \ \ \ \ \ \ \ \ \

Let $({\cal A},\Sigma)$ be a $C-$pair. We say that $({\cal A},\Sigma)$ is
complete if $\Sigma={\cal A}^{c}$.

{\bf Lemma 4.3. \ }Let $({\cal A},\Sigma )$ and $({\cal A}^{\prime },\Sigma
^{\prime })$ be $C-$pairs over $R.$ Assume that $\varphi :{\cal A}%
\rightarrow {\cal A}^{\prime }$ is a homomorphism of $R-$algebras and $(%
{\cal A},\Sigma )$ is complete. Then $\varphi :({\cal A},\Sigma )\rightarrow
({\cal A}^{\prime },\Sigma ^{\prime })$ is a $C-$homomorphism.

\bigskip {\bf Proof. \ }Let $\Lambda ^{\prime }\in \Sigma ^{\prime }.$ Then $%
\varphi ^{-1}(\Lambda ^{\prime })\in {\cal A}^{c}.$ Since $({\cal A},\Sigma
) $ is complete $\varphi ^{-1}(\Lambda ^{\prime })\in \Sigma .$ Therefore, $%
\varphi $ is a $C-$homomorphism.$\blacksquare $

We say that the $C-$pair $\left( {\cal A},\Sigma \right) $ is a $D-$pair if
for each $\Lambda \in \Sigma $ there exists $A\in D\left( {\cal A}\right) $
such that $\Lambda =\underline{A}.$

Let $\varphi :\left( {\cal A},\Sigma \right) \rightarrow \left( {\cal A}%
^{\prime },\Sigma ^{\prime }\right) $ be an injective $C-$homomorphism
between $D-$pairs. We say that $\varphi $ is a{\it \ domination, }or $\left(
{\cal A},\Sigma \right) $ is {\it dominated by }$\left( {\cal A}^{\prime
},\Sigma ^{\prime }\right) ${\it \ under }$\varphi $ or $\left( {\cal A}%
^{\prime },\Sigma ^{\prime }\right) $ {\it \ dominates }$\left( {\cal A}%
,\Sigma \right) $ under $\varphi $, if for each $\Lambda ^{\prime }\in {\cal %
A}^{\prime c},$ with $\varphi ^{-1}\left( \Lambda ^{\prime }\right) \in
\Sigma ,$ we have $\ \Lambda ^{\prime }\in \Sigma ^{\prime }$ and if $A\in
\left( \varphi ^{-1}\left( \Lambda ^{\prime }\right) \right) ^{c}$ , then $%
\varphi \left( A\right) \in \Lambda ^{\prime c}.$

{\bf Lemma 4.4. \ }Let \ $\varphi :\left( {\cal A},\Sigma \right)
\rightarrow \left( {\cal A}^{\prime },\Sigma ^{\prime }\right) $ be a
domination. Assume that $\left( {\cal A},\Sigma \right) $ is complete. Then,
the $C-$pair $\left( {\cal A}^{\prime },\Sigma ^{\prime }\right) $ is also
complete.

{\bf Proof. }\ Let {\bf \ }$\Lambda ^{\prime }\in {\cal A}^{\prime c}.$
Then, $\varphi ^{-1}\left( \Lambda ^{\prime }\right) \in {\cal A}^{c}.$
Since the $C-$pair $\left( {\cal A},\Sigma \right) $ is complete, $\varphi
^{-1}(\Lambda ^{\prime })\in \Sigma .$ As $\varphi $ is a domination, $%
\Lambda ^{\prime }\in \Sigma ^{\prime }.$ Therefore, $\left( {\cal A}%
^{\prime },\Sigma ^{\prime }\right) $ is complete.$\blacksquare $

{\bf Lemma\ 4.5. }\ Let $\varphi :({\cal A},\Sigma )\rightarrow ({\cal A}%
^{\prime },\Sigma ^{\prime })$\ be a domination. Assume that $A\in \overline{%
\Sigma }$ and $\varphi \left( A\right) =A^{\prime }\in \overline{\Sigma
^{\prime }}$. Then \ we have $\alpha _{1}^{\prime }\circ \varphi =\varphi
\circ \alpha _{1}$ and \ $\varphi \circ \underline{\alpha }_{2}=\underline{%
\alpha }_{2}^{\prime }\circ \varphi .$

{\bf Proof. }\ Let \ $y\in {\cal A}$. Then, $y=\alpha _{0}\left( y\right)
e+\alpha _{1}\left( y\right) +\underline{\alpha }_{2}\left( y\right) ,$
where $\alpha _{0}\left( y\right) \in R,$ $\alpha _{1}\left( y\right) \in A,$
and $\underline{\alpha }_{2}\left( y\right) \in \underline{A}^{2}.$ Thus, $%
\varphi (y)=\alpha _{0}\left( y\right) e^{\prime }+\varphi (\alpha
_{1}\left( y\right) )+\varphi (\underline{\alpha }_{2}\left( y\right) )$.
Since $\varphi \left( y\right) $ is an element of ${\cal A}^{\prime },$ we
have \ $\varphi \left( y\right) =\alpha _{0}^{\prime }\left( \varphi \left(
y\right) \right) e^{\prime }+\alpha _{1}^{\prime }\left( \varphi \left(
y\right) \right) +\underline{\alpha }_{2}^{\prime }\left( \varphi \left(
y\right) \right) ,$ where $\alpha _{0}^{\prime }\left( \varphi \left(
y\right) \right) \in R,$ $\alpha _{1}^{\prime }\left( \varphi \left(
y\right) \right) \in A^{\prime }$. and $\underline{\alpha }_{2}^{\prime
}\left( \varphi \left( y\right) \right) \in \underline{A^{\prime }}^{2}.$
Thus we have

\[
\alpha _{0}\left( y\right) e^{\prime }+\varphi (\alpha _{1}\left( y\right)
)+\varphi (\underline{\alpha }_{2}\left( y\right) )=\alpha _{0}^{\prime
}\left( \varphi \left( y\right) \right) e^{\prime }+\alpha _{1}^{\prime
}\left( \varphi \left( y\right) \right) +\underline{\alpha }_{2}^{\prime
}\left( \varphi \left( y\right) \right) .
\]%
Or

\[
\left( \alpha _{0}\left( y\right) -\alpha _{0}^{\prime }\left( \varphi
\left( y\right) \right) \right) e^{\prime }+\left( \varphi \circ \alpha
_{1}\left( y\right) -\alpha _{1}^{\prime }\circ \varphi \left( y\right)
\right) +\left( \varphi \circ \underline{\alpha }_{2}\left( y\right) -%
\underline{\alpha }_{2}^{\prime }\circ \varphi \left( y\right) \right) =0.
\]%
But $A^{\prime }\in D\left( {\cal A}^{\prime }\right) $ and $\varphi \left(
\underline{\alpha }_{2}\left( y\right) \right) \in \underline{A}^{\prime 2},$
Therefore:\ \ \ \ \ \ \ \ \ \ \ \ \ \ \ \ \ \ \ \ \ \ \ \ \ \ \ \ \ \ \ \ \
\ \ \ \ \ \ \ \ \ \ \ \ \ \ \ \ \ \ \ \ \ \ \ \ \ \ \ \ \ \ \
\[
\ \ \ \varphi \circ \alpha _{1}\left( y\right) =\underline{\alpha }%
_{1}^{\prime }\circ \varphi \left( y\right) \quad{ and }\quad
\varphi \circ \underline{\alpha }_{2}(y)=\underline{\alpha
^{\prime }}_{2}\circ \varphi (y).
\]%
Since $y\in {\cal A}$ is arbitrary, $\alpha _{1}^{\prime }\circ \varphi
=\varphi \circ \alpha _{1}$and $\varphi \circ \underline{\alpha }_{2}=%
\underline{\alpha ^{\prime }}_{2}\circ \varphi .\blacksquare $

{\bf Lemma\ 4.6. \ }Let\ $\varphi :\left( {\cal A},\Sigma \right)
\rightarrow \left( {\cal A}^{\prime },\Sigma ^{\prime }\right) $ be a
domination. Then, \ \underline{$\varphi $}:$\Sigma ^{\prime }\rightarrow
\Sigma $ is injective.

{\bf Proof. \ }Let $\Lambda ^{\prime },\Lambda ^{\prime \prime
}\in \Sigma ^{\prime },\Lambda \in \Sigma .$ Assume that $\varphi
^{-1}\left( \Lambda ^{\prime }\right) =\Lambda =\varphi
^{-1}\left( \Lambda ^{\prime \prime }\right) $. Then, \ since
$\varphi $ is a domination, for any $A\in \Lambda ^{\circ },$ we
have \ $\varphi \left( A\right) \in \Lambda ^{\prime
\circ }\cap \Lambda ^{\prime \prime \circ }.$ Thus, $\Lambda ^{\prime }=%
\underline{\varphi \left( A\right) }=\Lambda ^{\prime \prime }.$ Therefore,
\underline{$\varphi $} is injective.$\blacksquare $

\section{\protect\Large \ Smooth Structures}

Let $\left( {\cal A},\Sigma \right) $ be a $D-$pair. Assume that any 2
elements of $\overline{\Sigma }=\cup _{\Lambda \in \Sigma }\Lambda ^{\circ }$
are compatible with each other. Then $\left( {\cal A},\Sigma \right) $ is
called a {\it smooth pair}. Let $\Sigma $ be maximal with respect to the
above property. Then$\ \overline{\Sigma }$ is called a {\it smooth structure
on }${\cal A},$ and $\left( {\cal A};\overline{\Sigma }\right) $ is called a
{\it smooth algebra.}

\bigskip Let ${\cal A}$ be an $R-$algebra of \ type $C.$ Assume that ${\cal A%
}$ admits a sub-module $A$ of type $D$. By Zorn's lemma, there exists a
smooth structure $\overline{\Sigma }$ on ${\cal A}$ which contains $A.$

Assume that $\left( {\cal A},\Sigma \right) $ and $\left( {\cal A}^{\prime
},\Sigma ^{\prime }\right) $ are smooth $R-$pairs. A smooth morphism between
them is a morphism in the category $R-CP.$ It is clear that the class of all
smooth $R-$algebras (resp. $R-$pairs) together with smooth morphism between
them is a full subcategory of the category $R-CA$ (resp.$R-CP).$ This
category will be denoted by $R-SP.$

{\bf Lemma\ 5.1. \ }A necessary and sufficient condition for a smooth pair $%
\left( {\cal A};\Sigma \right) $ to be of polynomial type is that for every $%
\Lambda \in \Sigma $ there exists $A\in \Lambda ^{\circ }$ such that ${\cal %
A=\oplus }_{n=0}^{\infty }A^{n}.$

The proof is clear.

Let $\left( {\cal A};\Sigma \right) $ be a smooth pair. We say that $%
\overline{\Sigma }$ is a {\it complete smooth structure }on ${\cal A}$ if $%
\Sigma ={\cal A}^{c}.$ in this case $\left( {\cal A};\overline{\Sigma }%
\right) $ is called a {\it complete smooth algebra.}

{\bf Lemma\ 5.2.\ \ }Let{\bf \ }$\left( {\cal A};\overline{\Sigma }\right) $
be a complete smooth $R-$algebra and let $\left( {\cal A},\Sigma ^{\prime
}\right) $ be a smooth $R-$ pair. Then $\overline{\Sigma ^{\prime }}\subset
\overline{\Sigma }$.\

The proof is clear.

The ring $R$ is clearly an $R-$algebra. The singleton $\{\{0\}\}$ is the
unique smooth structure on $R.$ This structure will be denoted by $\left[ R%
\right] .$

Here we have the following simple lemma.

{\bf Lemma\ 5.3. \ }Assume that $\left( {\cal A};\overline{\Sigma }\right) $
is a smooth $R-$algebra. For $\Lambda \in {\cal A}^{c},$ we have $\Lambda
\in \Sigma $ if and only if $\pi _{\Lambda }^{0}:\left( {\cal A};\overline{%
\Sigma }\right) \rightarrow \left( R;\left[ R\right] \right) $ is smooth.

{\bf Proposition\ 5.4. }\ For{\bf \ }$n\geq 1$ the $%
\mathbb{R}
-$algebra ${\cal A}=C^{\infty }\left(
\mathbb{R}
^{{{ n}}}\right) ,$ admits a complete separated non-analytic
smooth structure which is the unique smooth structure on it.

{\bf Proof. \ }Let $x^{i}:%
\mathbb{R}
^{{{ n}}}\rightarrow
\mathbb{R}
$ denote the $i-$th projection and let $x:%
\mathbb{R}
^{{{n}}}\rightarrow
\mathbb{R}
^{{\it n}}$ denote the identity mapping and let $e$ denotes the constant map
$%
\mathbb{R}
^{{\it n}}:\rightarrow $ \{1\} $\subset
\mathbb{R}
.$ For $\lambda =\left( \lambda ^{1},\lambda ^{2},...,\lambda ^{{\it n}%
}\right) \in
\mathbb{R}
^{{\it n}}$ we set $M_{\lambda }=\Sigma _{k=1}^{n}%
\mathbb{R}
\cdot \left( {\it x}^{{\it k}}-\lambda ^{{\it k}}e\right) $ and$\ I_{\lambda
}=\Sigma _{k=1}^{n}{\cal A\cdot }\left( x^{{\it k}}-\lambda ^{{\it k}%
}e\right) .$ Clearly, $I_{\lambda }$ is generated by $M_{\lambda }.$ By a
well-known lemma whose statement will follow, one deduces that $I_{\lambda
}\in {\cal A}^{c}$. Assume that\ \ \ $\ \ \ \ \ \ \ \ \ \ \ \ \ \ \ \ \ \ \
\ \ \ \ \ \ \ \ \ \ \ \ \ \ \ \ \ \ \ \ \ \ \ \ \ \ \ \ \ \ \ \ \ \ \ \ \ \
\ \ \ \ \ \ \ \ \ \ \ \ \ \ \ \ \ \ \ \ \ \ \ \ \ \ \ \ \ $

\[
\ \ \Sigma _{\mid i\mid =k}\mu _{i}\left( x-\lambda e\right) ^{i}=\Sigma
_{\mid i\mid =k+1}\varphi _{i}\left( x-\lambda e\right) ^{i}\in M_{\lambda
}^{k}\cap I_{\lambda }^{k+1}
\]%
where $\mu _{i}\in
\mathbb{R}
$ and $\varphi _{i}\in {\cal A}$. Let $j=\left( j_{1},j_{2,}...,j_{n}\right)
\in
\mathbb{N}
^{{\it n}},$ $\mid ${\it j }$\mid =k$ and $D^{j}=\frac{\partial ^{k}}{%
\partial \left( x^{1}\right) ^{j_{1}}\partial \left( x^{2}\right)
^{j_{2}}...\partial \left( x^{n}\right) ^{j_{n}}}$. Then there exists a
non-zero constant $C_{j}\in
\mathbb{R}
$ such that\ \ \ \ \
\[
C_{j}\mu _{j}=D^{j}\left( \Sigma _{\mid i\mid =k}\mu _{i}\left( x-\lambda
e\right) ^{i}\right) \mid _{x=\lambda e}=D^{j}\left( \Sigma _{\mid i\mid
=k+1}\varphi _{i}\left( x-\lambda e\right) ^{i}\right) \mid _{x=\lambda
e}=0.
\]%
Therefore, $\mu _{j}=0$ and for each $k\in
\mathbb{N}
,$ $M_{\lambda }^{k}\cap {I_{\lambda }}^{k+1}=\{0\}.$ Hence,
$M_{\lambda { }}$ is a sub-module of type $D.$ Let $\lambda
,\lambda ^{\prime }\in
\mathbb{R}
^{{\it n}}.$ Then, clearly $M_{\lambda }$ and $M_{\lambda ^{\prime }}$ are
compatible with each other. Therefore, there exists a smooth structure on $%
{\cal A}$ which contains the set $\{M_{\lambda }\mid \lambda \in
\mathbb{R}
^{{\it n}}\}.$ This smooth structure will be denoted by $[%
\mathbb{R}
^{{\it n}}]$. Let $\Lambda \in {\cal A}^{c}$ and $\lambda =\left( \pi
_{\Lambda }^{0}\left( x^{1}\right) ,\pi _{\Lambda }^{0}\left( x^{2}\right)
,...,\pi _{\Lambda }^{0}\left( x^{n}\right) \right) .$ Then, $M_{\lambda
}\in \Lambda ^{\circ }$. Hence, each element of $[%
\mathbb{R}
^{{\it n}}]$ is equivalent to some $M_{\lambda }$ and the smooth algebra $%
\left( C^{\infty }\left(
\mathbb{R}
^{{\it n}}\right) ,[%
\mathbb{R}
^{{\it n}}]\right) $ is complete. Clearly, $\left( C^{\infty }\left(
\mathbb{R}
^{{\it n}}\right) ,[%
\mathbb{R}
^{{\it n}}]\right) $ is separated. Consider the function $\varphi :%
\mathbb{R}
^{{\it n}}\rightarrow
\mathbb{R}
$ defined by

\begin{equation*}
\varphi(t)=
\begin{cases}
e^{\frac{-1}{\sum_{1}^{n}(t^i)^2}}\quad if \quad\sum_{1}^{n}(t^i)^2\neq 0 & \\
0 \quad  \quad\quad\;\;     if\quad \sum_{1}^{n}(t^i)^2=0
\end{cases}
\end{equation*}

All the derivatives of $\varphi $ are zero at $\left(
0,0,...,0\right) .$ Therefore, $\left( C^{\infty }\left(
\mathbb{R}
^{{\it n}}\right) ,[%
\mathbb{R}
^{{\it n}}]\right) $ is not analytic. The uniqueness follows from
Lemma 5.2.$\blacksquare $

In a similar way without using the above lemma one can see that
the algebra of polynomial functions on an integral domain $R$ and
the algebra of entire functions on $\Bbbk^n (\Bbbk={\mathbb{R}}\;
or\; \Bbbk={\mathbb{C}})$ admit unique smooth structures. The
corresponding smooth algebras will be denoted by $\left( P\left(
R\right) ,[R]\right) $ and $\left( C^{\omega }\left( \Bbbk
^{n}\right) ,[\Bbbk ^{n}]\right) $, respectively. They are both
separated. The first is of polynomial type and the second is
analytic.

{\bf Lemma\ 5.5. \ }Let $f\in C^{\infty }\left(
\mathbb{R}
^{{\it n}}\right) $ and $\lambda \in
\mathbb{R}
^{{\it n}}.$ Then, there exist $n$ functions $g_{i}\in C^{\infty }\left(
\mathbb{R}
^{{\it n}}\right) ,$ $i=1,2,...,n,$ such that $f=f\left( \lambda \right)
e+\Sigma _{k=1}^{n}\left( x^{k}-\lambda ^{k}e\right) g_{k}$.

Now let $\Omega \subset
\mathbb{R}
^{{\it n}}$ be nonempty. Assume that $\lambda \notin \Omega .$ Define $%
\varphi _{\lambda }:%
\mathbb{R}
^{{\it n}}-\{\lambda \}\rightarrow
\mathbb{R}
,$ as follows%
\[
\varphi _{\lambda }=\Sigma _{k=1}^{n}\left( x^{k}-\lambda ^{k}e\right) ^{2}.
\]%
Let $W=\{\varphi _{\lambda }\mid \lambda \notin \Omega \}$ and let ${\cal B}=%
\frac{1}{W}{\cal A}$ be the localization of ${\cal A}=C^{\infty }(%
\mathbb{R}
^{{{ n}}}$ $)$ with respect to $W.$ Assume that $\Lambda \subset
\Lambda ^{\prime }$.\ Where, $\Lambda \in {\cal A}^{c}$ and
$\Lambda ^{\prime }\in {\cal B}^{c}.$Then, clearly $\lambda
=\left( \pi _{\Lambda }^{0}\left( x^{1}\right) ,\pi _{\Lambda
}^{0}\left( x^{2}\right) ,...,\pi
_{\Lambda }^{0}\left( x^{n}\right) \right) $ is an element of $%
\mathbb{R}
^{{\it n}}.$ and it is not difficult to prove that $\Lambda
^{\prime }$ is generated by $\Lambda $. Therefore, $\sigma :\Omega
\rightarrow {\cal B}^{c}$ is bijective and the $C-$pair (${\cal
B},\Omega )$ is complete. It is
separated if and only if $\Omega $ is dense in $%
\mathbb{R}
^{{\it n}}.$ It is called the {\it algebra of smooth functions associated
with} $\Omega .$Observe that the canonical injection $\iota :{\cal %
A\rightarrow B}$ is a domination.

{\bf Proposition\ 5.6. \ }Let $M$ be a closed sub-manifold of $%
\mathbb{R}
^{{\it n}}.$ Assume that ${\cal A=}C^{\infty }\left( M\right) $. Then, $%
{\cal A}^{c}=M$ and $C^{\infty }\left( M\right) $ is complete.

{\bf Proof. }\ Let{\bf \ }${\cal A}\left( M\right) $ denote the algebra of
smooth functions associated with $M$, \ and let $\varphi :{\cal A}\left(
M\right) \rightarrow {\cal A}$ be the restriction homomorphism. Since $(%
{\cal A}(M),M)$ is complete by Lemma 4.3. $\varphi :\left( {\cal
A}\left( M\right) ,M\right) \rightarrow \left( {\cal A}{,}{\cal
A}^{c}\right) $
is a $C-$homomorphism. Since $\varphi $ is surjective, by Lemma 4.2. iii, $%
\underline{\varphi }:{\cal A}^{c}\rightarrow \left( {\cal A}\left( M\right)
\right) ^{c}$ is injective. Clearly, \underline{$\varphi $}$\mid
_{M}=id_{M}. $ Therefore, ${\cal A}^{c}=M.\blacksquare $

From the above proposition and Whitney's embedding theorem we have the
following theorems.

{\bf Theorem\ 5.7.}\ \ Let{\bf \ }$M$ be a finite dimensional smooth
manifold. Then, $C^{\infty }\left( M\right) $ is complete.$\blacksquare $

{\bf Theorem\ 5.8.\ \ }Let $M$ and $N$ be two finite dimensional \ smooth
manifolds. Assume that $\varphi :C^{\infty }\left( N\right) \rightarrow
C^{\infty }\left( M\right) $ is a homomorphism of $%
\mathbb{R}
-$algebras. Then, there exists exactly one smooth mapping ${\small \Phi }%
:M\rightarrow N$ such that ${\small \Phi }^{\ast }=\varphi .$

{\bf Proof. \ }The uniqueness of ${\small \Phi }$ is clear. By the above
theorem $C^{\infty }\left( M\right) $ and $C^{\infty }\left( N\right) $ are
complete. Therefore, \underline{$\varphi $} is a map from $M$ into $N$.
Since for each $f\in C^{\infty }\left( N\right) $, ${\small \Phi }^{\ast
}\left( f\right) =\varphi \left( f\right) ,$ the mapping ${\small \Phi
:M\rightarrow N}$ is smooth.$\blacksquare $

We say that smooth $R$-pairs $\left( {\cal A},\Sigma \right) $, $\left(
{\cal A},\Sigma ^{\prime }\right) $ are consistent if each element of $%
\Sigma $ is compatible with each element of $\Sigma ^{\prime }$. In this
case $({\cal A},\Sigma \cup \Sigma ^{\prime })$ is also a smooth $R$-pair.

\section{ {\protect\Large Maximal Smooth Structure}}

{\bf Proposition\ 6.1. \ }Let $\left( I,\leq \right) $ be a directed set.
Assume that $\left( \left( {\cal A}_{i},\Sigma _{i}\right) ,\varphi
_{ji}\right) _{i,j\in I},i\leq j$ is a direct system of smooth pairs over $R$%
, where for each $i\leq j$, ${\cal A}_{i}$ is a subalgebra of ${\cal A}_{j}$
and $\varphi _{ji}$, the canonical injection of ${\cal A}_{i}$ into ${\cal A}%
_{j}$ is a domination. Then, the $R-$algebra ${\cal A=\cup }_{i\in I}{\cal A}%
_{i}$ is the underlying $R-$algebra of a unique smooth pair over $R$ such as
$\left( {\cal A},\Sigma \right) $ which satisfies the following conditions:

i) \ For each $i\in I$, $\varphi_{i}:\left( {\cal A}_{i},\Sigma _{i}\right)
\rightarrow\left( {\cal A},\Sigma\right) $ is a domination. (Here $%
\varphi_{i}:{\cal A}_{i}\rightarrow{\cal A}$ is the canonical injection.)

ii) \ If for some $p\in I$\ and for all $i\geq p,$ $\left( {\cal A}%
_{i},\Sigma _{i}\right) $ is separated (resp. analytic), then,
$\left( {\cal A},\Sigma \right) $ is separated (resp. analytic).

{\bf Proof. \ }Let $\mu $ be an element of $I$ and let $A\in \overline{%
\Sigma }_{\mu }.$ We are going to prove that $A$ is a sub-module of type $D$
for ${\cal A}$. Let $\Lambda =A\cdot {\cal A}$ and for $i\geq \mu ,$ let $%
\Lambda _{i}=A\cdot {\cal A}_{i}.$ Clearly, we have $\Lambda =\cup _{i\in
I}\Lambda _{i}.$ Since ${\cal A=\cup }_{i\in I}{\cal A}_{i},$ for each $y\in
{\cal A},$ there exists $i\in I$ such that $y\in {\cal A}_{i}.$ As $\Lambda
_{i}$ is an ideal of type $D$ for ${\cal A}_{i},$ $y=\pi _{\Lambda
_{i}}^{0}\left( y\right) +\pi _{\Lambda _{i}}\left( y\right) ,$ where $\pi
_{\Lambda _{i}}\left( y\right) \in \Lambda _{i}$. Since $\Lambda _{i}\subset
\Lambda ,$ $\pi _{\Lambda _{i}}\left( y\right) \in \Lambda $. Hence $y\in
R\cdot e+\Lambda .$ But $y$ is an arbitrary element of ${\cal A}$.
Therefore, ${\cal A=}R\cdot e+\Lambda .$ Now assume that for some $k\in
\mathbb{N}
,$\ \ \ \ \ \ \
\[
\Sigma _{\mid n\mid =k}\mu _{n}x^{n}=\Sigma _{\mid n\mid =k+1}y_{n}x^{n}\in
A^{k}\cap \Lambda ^{k+1},
\]%
\ \ \ \ where x$^{n}$=x$_{1}^{n_{1}}$x$_{2}^{n_{2}}$...x$_{q}^{n_{q}}$, \ n=$%
\Sigma _{l=1}^{q}$n$_{l}$, x$_{i}\in $A, $y_{n}\in {\cal A}$ and
$\mu
_{n}\in $R. Since ${\cal A=\cup }_{i\in I}{\cal A}_{i}$ and\ for i$\leq $j, $%
{\cal A}_{i}$\ is a subalgebra of ${\cal A}_{j},$ there exists $i\in I,$ $%
i\geq \mu ,$ such that all $y_{n}\in {\cal A}_{i}$. Thus in ${\cal A}_{i}$\
\ \ \ \ \ \ \ \ \ \ \ \ \ \ \ \ \ \ \ \ \ \ \ \ \ \ \ \ \ \ \
\[
\ \ \ \Sigma _{\mid n\mid =k}\mu _{n}x^{n}=\Sigma _{\mid n\mid =k+1{ }%
}y_{n}x^{n}\in A^{k}\cap \Lambda _{i}^{k+1}.
\]%
Since $A$ is a sub-module of type $D$ for ${\cal A}_{i}$ and generates the
ideal $\Lambda _{i}$\ \ \ \ \ \ \ \ \ \ \ \ \ \
\[
\ \ \ \ \ \ \ \ \Sigma _{\mid n\mid =k}\mu _{n}x^{n}=0=\Sigma _{\mid n\mid
=k+1}y_{n}x^{n}.
\]%
Therefore, $A$ is a sub-module of type $D$ for ${\cal A}$. On the other
hand, let $A$ and $A^{\prime }$ be two \ equivalent sub-modules of type $D$
for $\left( {\cal A}_{\mu },\Sigma _{\mu }\right) $. Since each $\varphi
_{j\mu }:{\cal A}_{\mu }\rightarrow {\cal A}_{j}$ is a domination, $\
A^{\prime }\cdot {\cal A}_{j}=A\cdot {\cal A}_{j}=\Lambda _{j}.$ Hence, $%
A^{\prime }\cdot {\cal A=}\cup _{i\in I}\left( A^{\prime }\cdot {\cal A}%
_{i}\right) =\cup _{i\geq \mu }\Lambda _{i}=\Lambda .$ Therefore, $A$ and $%
A^{\prime }$ are equivalent in ${\cal A}$. Let \ $C\in \overline{\Sigma }%
_{\mu }$ be strongly compatible with $\ B\in \overline{\Sigma }_{\mu }.$ By
definition $\pi _{\underline{B}}(C)\thickapprox B$ and $\pi _{\underline{C}%
}(B)\thickapprox C$ in ${\cal A}_{\mu }.$ By the above, $\pi _{\underline{B}%
}(C)\thickapprox B$ and $\pi _{\underline{C}}(B)\thickapprox C$ in ${\cal A}%
. $ Hence, if $D$ and $E$ are compatible in ${\cal A}_{\mu }$ they are
compatible in ${\cal A}.$ Therefore, there exists a smooth pair $\left(
{\cal A},\Sigma \right) $ where \ $\Sigma =\cup _{\Lambda \in \Sigma
_{i}}\varphi \left( \Lambda \right) $ for all $i\in I.$ Clearly, each $%
\varphi _{i}:({\cal A}_{i},\Sigma _{i})\rightarrow ({\cal A},\Sigma )$ is a
domination. The uniqueness of this smooth pair is trivial.

Assume that for all $i\geq \nu ,\left( {\cal A}_{i},\Sigma _{i}\right) $ is
separated. Let $y\in \cap _{\Lambda \in \Sigma }\Lambda $. Then, there
exists $i\geq \nu ,$ such that $y\in {\cal A}_{i}$. Thus $y=\varphi
_{i}{}^{-1}\left( \cap _{\Lambda \in \Sigma }\Lambda \right) =\cap _{\Lambda
\in \Sigma }\varphi _{i}^{-1}\left( \Lambda \right) =\cap _{\Lambda \in
\Sigma _{i}}\Lambda =\{0\}.$ Therefore, $\left( {\cal A},\Sigma \right) $ is
separated.

In the same way we see that if for all $i\geq\nu,\left( {\cal A}%
_{i},\Sigma_{i}\right) $ is analytic, then $\left( {\cal A}%
,\Sigma\right) $ is analytic.$\blacksquare$

Let $\Gamma $ be a subcategory of the category $R-SP.$ An object $\left(
{\cal A},\Sigma \right) $ is called {\it maximal in }$\Gamma $ if every
domination $\varphi :({\cal A},\Sigma )\rightarrow \left( {\cal A}^{\prime
},\Sigma ^{\prime }\right) $ which is in $\Gamma $ is an isomorphism.

{\bf Proposition\ 6.2. \ }The smooth $%
\mathbb{R}
-$algebra $\left( C^{\infty }\left(
\mathbb{R}
^{{\it n}}\right) ,\left[
\mathbb{R}
^{{\it n}}\right] \right) $ is not separatedly maximal (is not maximal in
the category of separated smooth $%
\mathbb{R}
-$algebras.)

{\bf Sketch of the Proof. \ }For simplicity we assume that $n=1.$ Let $%
\lambda \in
\mathbb{Q}
$ and $\mu \in
\mathbb{R}
-%
\mathbb{Q}
.$ Define $\varphi _{\lambda },\psi _{\mu }:%
\mathbb{R}
\rightarrow
\mathbb{R}
$ as follows
\begin{equation*}
\varphi_{\lambda}(t)=
\begin{cases}
{\frac{1}{t-\lambda}}\quad if \quad t\in\mathbb{R}-\mathbb{Q}& \\
0 \quad  \quad\quad\;\;     if\quad t\in\mathbb{Q}
\end{cases}
\end{equation*}

\begin{equation*}
\psi_{\mu}(t)=
\begin{cases}
{\frac{1}{t-\mu}}\quad if \quad t\in\mathbb{Q}& \\
0 \quad  \quad\quad\;\;     if\quad t\in\mathbb{R}-\mathbb{Q}
\end{cases}
\end{equation*}

Let ${\cal A}$ be the subalgebra of $%
\mathbb{R}
^{%
\mathbb{R}
}$ generated by $C^{\infty }\left(
\mathbb{R}
\right) \cup \left\{ \varphi _{\lambda }\mid \lambda \in
\mathbb{Q}
\right\} \cup \left\{ \psi _{\mu }\mid \mu \notin
\mathbb{Q}
\right\} $. As we have done in Proposition 5.4 one can check easily that $%
{\cal A}$ admits a unique complete separated smooth structure $\overline{%
\Sigma }$ such that the canonical injection $\left( C^{\infty }\left(
\mathbb{R}
\right) ,%
\mathbb{R}
\right) \rightarrow \left( {\cal A},\Sigma \right) $ is a domination which
is not an isomorphism.

The proof in the general case is the same.$\blacksquare$

{\bf Theorem\ 6.3. \ }Every separated smooth{\bf \ }$R-$pair ($R-$algebra)
is dominated by a maximal one.

{\bf Proof. \ }Let $\left( {\cal A},\Sigma \right) $ be a separated smooth $%
R-$pair. Assume that $\Omega $ is the set consisting of all separated smooth
$R-$pairs which dominate $\left( {\cal A},\Sigma \right) $. Since $\left(
{\cal A},\Sigma \right) \in \Omega ,$ $\Omega $ is not empty. Now we order $%
\Omega $ by domination. By Proposition 6.1, each chain in $\Omega $ has an
upper bound in $\Omega .$ By Zorn's lemma $\Omega $ has a maximal element $%
\left( {\cal A}^{\prime },\Sigma ^{\prime }\right) $ , which dominates $%
\left( {\cal A},\Sigma \right) .\blacksquare $ \ \ \ \ \ \ \ \ \ \ \ \ \ \ \
\ \ \ \ \ \ \ \

\section{\protect\Large \ Smooth Structures on Rings}

Let $R$ be a commutative ring with identity 1. We say that $R$ is a {\it %
semi-integral domain}, if for $x,y\in R,$ and $n\in
\mathbb{N}
,$ the relations $x\neq 0$ and $x^{n}\left( xy+1\right) =0$ imply
that $x$ is invertible. $R$ is called {\it proper} if there is a
non-unit $t\in R$ such that $1+t$ is also a non-unit. Clearly,
each integral domain is a semi-integral domain. Furthermore, we
have the following

{\bf Lemma 7.1.}\ \ Let $R$ be a commutative ring with identity 1. Then

1) $R$ is a semi-integral domain if and only if it contains neither
nilpotent nor idempotent elements.

2) $R$ is proper if and only if it is not local.

3) Let $R$ be a semi-integral domain and $R^{\prime }\subset R$ a
sub-ring. Then $R^{\prime }$ is also a semi-integral domain.

4) Any finite semi-integral domain is a field.

\bigskip {\bf Proof.} 1) It follows from the definition that if $R$ has a
nilpotent or an idempotent element, then it cannot be a semi-integral
domain. Now assume that it is not a semi-integral domain. Then there exist $%
x,y\in R,$ and $n\in
\mathbb{N}
$ , such that $x\neq 0$ and $x^{n}\left( xy+1\right) =0.$ But $x$ is not a
unit. If $x^{n}=0$, there is nothing to prove. Otherwise let $z=-xy.$ So $%
z^{n}(z^{n}-1)=z^{n}(z-1)(z^{n-1}+z^{n-2}+...+z+1)=0$. Since $z$ is not a
unit $\ $it is a non-trivial nilpotent or $\ z^{n}$ is an idempotent.

2) Let $(R,{\frak m})$ be local. Then the non-units are precisely the
elements of ${\frak m}$. Clearly for $x\in {\frak m}$ , $1+x$ is a unit.
Thus $R$ is not proper.

Conversely, let ${\frak m}_{1}$ and ${\frak m}_{2}$ be two maximal ideals
and let ${\frak m}_{1}\neq {\frak m}_{2}.$ Then, since the ideal generated
by ${\frak m}_{1}\cup {\frak m}_{2}$ is the ring $R$, there exist $x\in
{\frak m}_{1}$ and $y\in {\frak m}_{2}$ such that $x-y=1.$ Therefore $1+y=x.$

3) The proofs are clear.

4) Let $R$ be a finite semi-integral domain with $n$
elements,and let $0\neq a\in R.$ Then for some $0\leq i\leq n$ we have $%
a^{n+1}=a^{i}$. Therefore, $a^{i}(1-a^{n-i+1})=0.$ Or
$a^{i}(1+a(-a^{n-i}))=0.$ Since $R$ is a semi-integral domain $a$
is invertible.$\blacksquare $

{\bf Lemma 7.2. \ }Let $X$ be a non-coarse connected topological
space.
Then, $R=C\left( X,%
\mathbb{R}
\right) $ is a proper semi-integral domain.

{\bf Proof.} \ Let $f,g\in R$ , $f\neq 0$ and let $Y=\{x\in X\mid f\left(
x\right) =0\}.$ Since $f\neq 0,$ $Y\neq X.$ Assume that $Y\neq \phi .$ Then,
there exists a point $x_{0}$ $\in Y,$ such that $Y$ \ is not a neighbourhood
of it. Let $U$ be an open set containing $x_{0}$ such that for all $x\in U,$
$f\left( x\right) g\left( x\right) +1\geq \frac{1}{2}.$ Since $%
U\varsubsetneq Y,$ there exists a point $x\in U$ such that $f\left( x\right)
\neq 0.$ Therefore, $f\left( x\right) ^{n}\left( f\left( x\right) g\left(
x\right) +1\right) \neq 0$ for all $n\in
\mathbb{N}
.$ It is clear that $R$ is proper.$\blacksquare $

Let $R$ be a commutative ring with identity 1. A subset $\Omega $ of $R$ is
called {\it absorbing }(resp.{\it \ strongly absorbing }){\it with respect to%
} $\lambda \in \Omega ,$ if for each $\lambda \neq $ $\alpha $ in
$R$ there exists $\beta \in R$\ (resp. $\beta \in \Omega $) such
that for $\alpha \neq 0,\ $we have $\alpha \beta \neq 0$ and
$\lambda +\beta \alpha \in \Omega .$ A non-empty subset $\Omega $
of $R$ is called{\it \ absorbing }(resp.{\it strongly absorbing})
if it is absorbing (resp.strongly absorbing) with respect to all
of its elements. An absorbing subset $\Omega \subset R$ is
called{\it \ proper} if for each invertible element $\lambda $ in
$\Omega $ there exists a non-invertible element $x$ of $R$ such
that $\lambda +x$ is a non-invertible element of $\Omega .$

Observe that all ideals are absorbing subsets and the
intersection of two absorbing subsets of a ring may be void.

In the following $x:\Omega \rightarrow R$ denotes the inclusion map and $%
e:\Omega \longrightarrow R$ is the constant mapping $t:\longmapsto 1$.

{\bf Lemma 7.3. \ }Let{\bf \ }$R$ be a semi-integral domain which
is not a field, and let $\Omega $ be an absorbing subset of $R.$
Assume that $\lambda \in \Omega ,$ $0\neq \eta \in R$ . Let
$\varphi _{\lambda ,\eta }:\Omega \rightarrow R$ be defined as
follows
\begin{equation*}
\varphi_{\lambda,\eta}(t)=
\begin{cases}
\eta \quad if \quad t=\lambda& \\
0 \quad  \quad\quad\;\;     if\quad t\neq\lambda
\end{cases}
\end{equation*}
Then there exists no function $\psi :\Omega \rightarrow R$, such that $%
\varphi _{\lambda ,\eta }$ is written as

\begin{center}
\[
(\ast )\ \ \ \ \ \ \ \ \ \ \ \ \ \ \ \ \ \ \ \ \ \varphi _{\lambda ,\eta
}=\eta e+\left( x-\lambda e\right) \psi
\]
\end{center}

{\bf Proof. \ }Assume that there exists $\psi :\Omega \rightarrow
R$ \ such that $\varphi _{\lambda ,\eta }$ can be written as
above. Since $\Omega $ is absorbing there exists $\beta \in R,$
such that
$\lambda \neq \lambda +\eta ^{2}\beta \in \Omega .$
Hence$\ \ \ \ \ \ \ $

\[
\ 0=\varphi _{\lambda ,\eta }\left( \lambda +\eta ^{2}\beta \right) =\eta
+\eta ^{2}\beta \cdot \psi \left( \lambda +\eta ^{2}\beta \right) =\eta
\lbrack 1+\eta \beta \cdot \psi \left( \lambda +\eta ^{2}\beta \right) ].
\]%
Since $R$ is a semi-integral domain, the equality $\eta \lbrack 1+\eta \beta
\cdot \psi \left( \lambda +\eta ^{2}\beta \right) ]=0$ and $\eta \neq 0$
imply that $\eta \left[ -\beta \psi \left( \lambda +\eta ^{2}\beta \right) %
\right] =1.$ Thus, $\eta $ is a unit. Now the definition of $\varphi
_{\lambda ,\eta }$ implies that for each $\alpha \in R$ and $\alpha \neq 0,$
there exists some $\beta \in R$ such that $\lambda +\alpha \beta \in \Omega
, $ and $\alpha \beta \cdot \psi \left( \lambda +\alpha \beta \right) =-\eta
.$ Since $R$ is not a field, this is a contradiction.$\blacksquare $

{\bf Lemma 7.4.\ \ }Let $R$ be a semi-integral domain and let
$\Omega $ be an absorbing subset of $R.$ Then if for $k\in
\mathbb{N}
,$ $0\neq \mu \in R,$ $\lambda \in \Omega $ and $g\in R^{\Omega }$ the
function $h:\Omega \longrightarrow R$ defined by

\[
h(x)=(x-\lambda e)^{k}(\mu e+(x-\lambda e)g(x))
\]%
is identically zero, then $R$ is a field.

{\bf Proof.\ \ }Assume that $h$ is identically zero. Since{\bf \ }$\Omega $
is absorbing there exists $\eta \in R$ such that $\lambda \neq x=\lambda
+\mu ^{2}\eta \in \Omega .$ So $\mu ^{2k}\eta ^{k}(\mu +\eta \mu
^{2}g(\lambda +\mu ^{2}\eta ))=0.$ Therefore

\[
(\mu \eta )^{2k+1}(1+\mu \eta g(\lambda +\mu ^{2}\eta ))=0.
\]%
Since $R$ is a semi-integral domain and $\mu \eta \neq 0,$ $\mu $ is a unit.
Without any loss of generality we assume that $\mu =1.$ Then the relation
\[
{\ }h(x)=(x-\lambda e)^{k}(e+(x-\lambda e)g(x))
\]%
and the fact that $R$ is a semi-integral domain imply that $R$ is a field.

{\bf Definition 7.5. \ }Let $R$ be a semi-integral domain which is
not a field and let $\Omega $ be an absorbing subset of $R$. A
{\it smooth structure }on{\it \ }$\Omega $ is a smooth pair
$\left( {\cal A},\Sigma \right) $ with the following properties:

i)\ The algebra ${\cal A}$ is a sub-algebra of $R^{\Omega }$ and $\Sigma
=\sigma \left( \Omega \right) .$ In the following we identify $\Sigma $ and $%
\Omega $ by $\sigma .$

ii) The inclusion map $x:\Omega \longrightarrow R$ is in ${\cal A}.$

iii) \ The smooth pair $\left( {\cal A}{,}\Sigma \right) $ is
separated.

iv) The pair\ $\left( {\cal A}{,}\Sigma \right) $\ is separatedly
maximal.

{\bf Theorem\ 7.6. \ }Let $R$ and $\Omega $ be as above. Then,

1) There exists a unique smooth structure $\left( {\cal A},\Sigma \right) $
on $\Omega .$

2) If $\Omega ^{\prime }$ is any other absorbing subset of $R$ contained in $%
\Omega ,$ and if $\left( {\cal A}^{\prime },\Sigma ^{\prime }\right) $ is
the smooth structure on $\Omega ^{\prime },$ then the restriction of each
element of ${\cal A}$ to $\Omega ^{\prime }$ is an element of ${\cal A}%
^{\prime }$ If , $R$ is an integral domain, or for some $\lambda \in \Omega
^{\prime },$ $\cap _{n=0}^{\infty }I_{\lambda }^{n}=\{0\}.$

{\bf Proof.} \ Let ${\cal P}\left( R\right) $ denote the $R-$algebra of
polynomial functions on $R.$ Then, $\left( {\cal P}\left( R\right) ,\Omega
\right) $ is a smooth pair which satisfies conditions i-iii above. It is
clear that each smooth structure on $\Omega $ must dominate $\left( {\cal P}%
\left( R\right) ,\Omega \right) .$ By Theorem 6.3 there exists a separatedly
maximal smooth pair $\left( {\cal A},\Omega \right) $ which dominates $%
\left( {\cal P}\left( R\right) ,\Omega \right) .$ Clearly, $\left( {\cal A}%
,\Omega \right) $ satisfies all the conditions i-iv above. We are going to
prove that $\left( {\cal A},\Omega \right) $ is the unique smooth pair which
satisfies all the above conditions. Suppose that $\left( {\cal B},\Omega
\right) $ is another smooth structure on $\Omega .$ Let ${\cal A\vee B}$
denote the sub-algebra of $R^{\Omega }$ generated by ${\cal A\cup B}$.
Consider the $C-$pair $\left( {\cal A\vee B},\Omega \right) .$ For $\lambda
\in \Omega ,$ let $M_{\lambda }=R\cdot \left( x-\lambda e\right) $ and $%
I_{\lambda }=\left( x-\lambda e\right) \left( {\cal A\vee B}\right) .$
Assume that $y\in {\cal A\vee B}$ . Then, there exists $x_{i}\in {\cal A}$, $%
y_{i}\in {\cal B}$, $i=1,2,...,n,$ such that $\ \ \ \ \ \ \ \ \ \ \ \ \ \ \
\ \ $%
\begin{eqnarray*}
y &=&\Sigma _{i=1}^{n}x_{i}.y_{i}=\Sigma _{i=1}^{n}(x_{i}(\lambda
)e+(x-\lambda e)\overline{x_{i}})(y_{i}(\lambda )e+(x-\lambda e)\overline{%
y_{i}}) \\
&=&\Sigma _{k=1}^{n}x_{k}(\lambda )y_{k}(\lambda )e+(x-\lambda e)\Sigma
_{i=1}^{n}(x_{i}(\lambda )\overline{y_{i}}+y_{i}(\lambda )\overline{x_{i}}%
+(x-\lambda e)\overline{x_{i}}.\overline{y_{i}}) \\
&\in &R.e+I_{\lambda }=R.e+M_{\lambda }.({\cal A\vee B}).
\end{eqnarray*}%
Therefore, \ for each $\lambda \in \Omega ,$ $M_{\lambda }$\ $\in C\left(
{\cal A\vee B}\right) $. Now assume that for some $\lambda ,\mu \in R,$ $%
0\neq k\in
\mathbb{N}
,$ and some $g\in {\cal A\vee B}$, we have $\mu (x-\lambda
e)^{k}=(x-\lambda e)^{k+1}g$. Then, $h(x)=\left( x-\lambda
e\right) ^{k}$\ $\left( \mu e+\left( x-\lambda e\right)
(-g\right) )=0$. By Lemma 7.6, $\mu \neq 0$ implies that the above
equality is impossible. Thus, for each $\lambda \in \Omega ,$
$M_{\lambda }^{k}\cap I_{\lambda }^{k+1}=0.$ Therefore, $\left(
{\cal A\vee B},\Omega \right) $ is a smooth pair\ which
dominates\ $\left( {\cal A},\Sigma \right) .$ Since\ $\left(
{\cal A},\Sigma \right) $\ is maximal\ \ we have ${\cal A=B}$.

Now assume that $\Omega ^{\prime }$ is an absorbing\ subset of
$R$ included in $\Omega .$ Let $y$ and $z$ be elements
of ${\cal A}$ and let\ $\lambda \in \Omega ^{\prime }.$ Then there exists $%
\overline{y}$,$\overline{z}\in {\cal A}$, such that $y=y\left(
\lambda \right) e+\left( x-\lambda e\right) \overline{y}$\ and
$z=z\left( \lambda \right) e+\left( x-\lambda e\right)
\overline{z}.$ Assume that $y\mid _{\Omega ^{\prime }}=$ \ $z\mid
\Omega ^{\prime }.$ Then, the above equalities imply that $\left(
x-\lambda e\right) \left( \overline{y}\mid _{\Omega ^{\prime
}}-\overline{z}\mid _{\Omega ^{\prime }}\right) $\ $=0.$
Assume that $R$ is an integral domain Then%
\begin{equation*}
\overline{z}\ \left( t\right) -\overline{y}\ \left( t\right)=
\begin{cases}
\overline{z}\ \left( \lambda\right) -\overline{y}\ \left(\lambda \right) \quad if \quad t=\lambda \\
0 \quad  \quad\quad\;\;     if\quad t\neq \lambda
\end{cases}
\end{equation*}
For $t\in \Omega ^{\prime }$.

\ \ \ \ But $\Omega ^{\prime }$ is an absorbing subset of $R.$ Thus Lemma
7.3 implies that $\overline{y}$\ $\mid _{\Omega ^{\prime }}=$\ $\overline{z}$%
\ $\mid _{\Omega ^{\prime }}.$\ Now assume that $\cap _{n=1}^{\infty
}I_{\lambda }^{n}=0.$ Then if $\overline{z}\ \left( t\right) -\overline{y}\
\left( t\right) \in $\ $\cap _{n=1}^{\infty }I_{\lambda }^{n}$, $\overline{z}%
\ \left( t\right) =\overline{y}\ \left( t\right) .$ Otherwise, by Lemma 7.5 $%
\overline{z}\ \left( t\right) =\overline{y}\ \left( t\right) .$ Therefore,
the restriction of elements of ${\cal A\ }$to $\Omega ^{\prime }$ is the
underlying $R-$algebra of a separated smooth pair which admits $M_{\lambda }$%
\ ,\ $\lambda \in \Omega ^{\prime }$\ as modules of type $D.$ By the unicity
of smooth structure on absorbing subsets of $R,$\ this smooth algebra is
included in $\left( {\cal A},\Sigma \right) $\ .\ $\blacksquare $

Let $\Omega $ be as above. Each $y\in {\cal A}$ is called a{\it \ smooth
function} on $\Omega $. Let $\lambda \in \Omega ,$ $k\in
\mathbb{N}
.$ Then $y\in {\cal A}$ can be written uniquely as

 $y=y\left( \lambda \right) e+a_{1}\left( x-\lambda e\right)
+a_{2}\left( x-\lambda e\right) ^{2}+...+a_{k}$\ $\left(
x-\lambda e\right) ^{k}+z\cdot
\left( x-\lambda e\right) ^{k+1},$where $a_{i}$\ $\in R$ and $z\in {\cal A}$%
. The element $k!\cdot a_{k}$\ of $R$ is called the ${\it k-th}$ {\it %
derivative of }$y${\it \ at }$\lambda $ and is denoted by\ $\frac{d^{n}y}{%
dx^{n}}\left( \lambda \right) $. The function $\frac{d^{k}y}{dx^{k}}:\Omega
\rightarrow R$ defined by $\frac{d^{k}y}{dx^{k}}$\ :$\lambda \longmapsto
\frac{d^{k}y}{dx^{k}}\left( \lambda \right) $\ is called the $k-th$ {\it %
derivative }of $y.$ Clearly we have $\frac{d^{k}\left( y\cdot z\right) }{%
dx^{k}}=\Sigma _{n=0}^{k}$\ $C_{k}^{n}$\
$\frac{d^{k-n}y}{dx^{k-n}}\cdot \frac{d^{n}z}{dx^{n}}$\ .

{\bf Important}\ {\bf Remark}\ {\bf 7.7. \ }Let $\ \Omega $ and
${\cal A}$ be as above. Suppose that $S\subset R^{\Omega }$\ .
The above theorem and its proof show that to prove that \ $S$ is
included in ${\cal A}$, it is sufficient to construct an
$R-$algebra\ \ ${\cal A}^{\prime }\subset R^{\Omega }$\
containing $S\cup \left\{ x\right\} $\ and prove that for each
$y\in {\cal A}^{\prime }$ and each $\lambda \in \Omega ,$ there exists $%
y_{\lambda }$\ $\in {\cal A}^{\prime }$ such that $y=y\left( \lambda \right)
e+\left( x-\lambda e\right) y_{\lambda }.$\

{\bf Theorem 7.8. \ }Let $\Omega $ be an absorbing subset of $R$.
Assume
that $\varphi :%
\mathbb{N}
\rightarrow \Omega ${\it \ } is bijective. Suppose that $f_{i}$\ :$%
\mathbb{N}
\rightarrow
\mathbb{N}
,${\it \ i = 1, 2, ...,} is a sequence of unbounded increasing functions.
Let $\left( a_{i}\right) _{i\in
\mathbb{N}
}$\ be a sequence of elements of $R$. Then the function $H:\Omega
\rightarrow R$ given by\ \ \ \ \ \ \ $\ .$

\[
\ H=\Sigma _{i=0}^{\infty }[a_{i}\Pi _{k=0}^{i}(x-\varphi (k)e)^{f_{k}(i)}],
\]%
is a smooth function on $\Omega $.\

{\bf Proof.} \ Let ${\cal A}$ be the sub-algebra of $R^{\Omega }$\ generated
by all\ functions of the above form. Assume that $\lambda \in \Omega $ and $%
H=\Sigma _{i=0}^{\infty }[a_{i}\Pi _{k=0}^{i}\left( x-\varphi \left(
k\right) e\right) ^{f_{k}\left( i\right) }]$ is an element of ${\cal A}$.
Then, there exists\ $n\in
\mathbb{N}
$ such that $\lambda =\varphi \left( n\right) $\ and

\[
H=\Sigma _{i=0}^{h-1}\left[ a_{i}\Pi _{k=0}^{i}\left( x-\varphi \left(
k\right) e\right) ^{f_{k}\left( i\right) }\right] +\left( x-\lambda e\right)
\Sigma _{i=h}^{\infty }\left[ a_{i}\Pi _{k=0}^{i}\left( x-\varphi \left(
k\right) e\right) ^{\overline{f}_{k}\left( i\right) }\right]
\]%
where

\begin{equation*}
\overline{f_{k}}(i)=
\begin{cases}
f_{k}(i)\quad if \quad k\neq n \\
f_{n}(i)-1 \quad  \quad\quad\;\;     if\quad k=n
\end{cases}
\end{equation*}
and h is such that $f_{n}\left( h\right) \neq 0.$ Since the first part of
the above sum is a polynomial\ in $x,$ for some $\mu \in R$ and some $%
P\left( x\right) \in {\cal P}\left( R\right) $\ it can be written as

\[
\Sigma _{i=0}^{h-1}a_{i}\Pi _{k=0}^{i}\left( x-\varphi \left( k\right)
e\right) ^{f_{k}\left( i\right) }=\mu e+\left( x-\lambda e\right) P\left(
x\right)
\]%
Therefore,

\[
H=\mu e+\left( x-\lambda e\right) [P\left( x\right) +\Sigma _{i=h}^{\infty
}[a_{i}\Pi _{k=0}^{i}\left( x-\varphi \left( k\right) e\right) ^{f_{k}\left(
i\right) }].
\]%
By the above remark $H$ is a smooth function.$\blacksquare $

{\bf Remark 7.9. \ }Assume that the characteristic of $R$ is
$k\neq 0.$ Then clearly all derivatives of the non-constant
function $x^{n}$ are identically zero. It is clear that the
function $\varphi =\frac{(x-e)x}{2}$ is not smooth. But the
function $2\varphi $ is smooth.

From now on the $R$-algebra of smooth functions on $\Omega \subset R$ will
be denoted by $C^{\infty }(\Omega )$.

{\bf Lemma 7.10.}\ \ Let $R$ be a proper semi-integral domain and
let $f\in
C^{\infty }(R)$ be a non-zero non-invertible function. Then there exists $%
x\in R$ such that $f(x)\neq 0$ and is non-invertible.

{\bf Proof.\ }\ The proof is by absurd. Assume that for each $x\in R$, $f(x)$
is zero or invertible. Then

1) Let $f(0)$\ $=0$. Then $f(x)=xg(x).$ Hence if $x$ is not unit $f(x)$ must
be zero. Let $f(\lambda )$ be invertible. Then $f(x)=\mu +(x-\lambda )h(x).$
Where $\mu $ is a unit. Assume that $x\in R$ is a non-unit such that $1+x$
is also non-unit. Then $-\lambda x$ and $-\lambda x-\lambda $ are also
non-unit. Since $f(-\lambda x)=0$ , $\mu +(-\lambda )(x+1)h(x)=0.$ Which is
absurd.

2) Let $f(0)=\alpha $ be a unit. Then $f(x)=\alpha +xg(x).$ Therefore if $x$
is not unit $f(x)$ is invertible. Assume that $\lambda $ is such that $%
f(\lambda )$ is zero. Then $\lambda $ is invertible and $f(x)=(x-\lambda
)k(x).$ Now assume that $x$ and $1+x$ are not invertible. Then $\lambda x$
is also non-invertible. But then $f(-\lambda x)=(-\lambda )(x+1)k(-\lambda
x) $ is invertible. Which is absurd.$\blacksquare $

{\bf Proposition\ 7.11.}\ \ Let $R$ be a proper semi-integral domain. Then $%
C^{\infty }(R)$ is also a proper semi-integral domain.

The proof is an immediate consequence of the above lemma.\ \ \ \ \ \ \ \ \ \
\ \ \ \ \ \ \ \ \ \ \ \ \ \ \ \ \ \ \ \ \ \ \ \ \ \ \ \ \ \ \ \ \ \ \ \ \ \
\ \ \ \ \ \ \ \ \ \ \ \ \ \ \ \ \ \ \ \ \ \ \ \ \ \ \ \ \ \ \ \ \ \ \ \ \ \
\ \ \ \ \ \ \ \ \ \ \ \ \ \ \ \ \ \ \ \ \ \ \ \ \ \ \ \ \ \ \ \ \ \ \ \ \ \
\ \ \ \ \ \ \ \ \ \ \ \ \ \ \

Let \ $R$ be a local integral domain which is not a field. Suppose that $M$
is the maximal ideal of $R.$ Assume that $\psi :R\rightarrow R$ is a smooth
function. Then, clearly for each $\lambda \in M$ , the function $\frac{1}{%
1+\lambda \psi }$ is a smooth function.\ \ \ \ \ \ \ \ \ \ \ \ \ \ \ \ \ \ \
\ \ \ \ \ \ \ \ \ \ \ \ \ \ \ \ \ \ \ \ \ \ \ \ \ \ \ \ \ \ \ \ \ \ \ \ \ \
\ \ \ \ \ \ \ \ \ \ \ \ \ \ \ \ \ \ \ \ \ \ \ \ \ \ \ \ \ \ \ \ \ \ \ \ \ \
\ \ \ \ \ \ \ \ \ \ \ \ \ \ \ \ \ \ \ \ \ \ \ \ \ \ \ \ \ \ \ \ \ \ \ \ \ \
\ \ \ \

{\bf Proposition\ 7.12. }\ Let $R$ be an semi-integral domain
which is not a field. Assume that $\Omega \subset R$ is an
absorbing subset and let $y$ be a smooth function on $\Omega $.
Then

1) For all $n\geq1,\frac{d^{n}y}{dx^{n}}:\Omega\rightarrow R$ is smooth.

2) If $\frac{1}{y}$ is defined on an absorbing subset $\Omega ^{\prime
}\subset $ $\Omega $, then as a function from $\Omega ^{\prime }$into$R$ is
also smooth and $\frac{d\left( \frac{1}{y}\right) }{dx}=-\frac{1}{y^{2}}%
\cdot \frac{dy}{dx}$ .

3) If $y\left( \Omega \right) $ is included in an absorbing subset of $R$,
say $\Omega ^{\prime }$ and $z\in {\cal A}_{\Omega ^{\prime }},$ then $%
z\circ y$ is a smooth function and $\frac{d\left( z\circ y\right) }{dx}=%
\frac{dz}{dx}\circ y\times \frac{dy}{dx}.$

{\bf Proof.\ }1) Let ${\cal A}_{\Omega }$ be the $R-$algebra generated by
all th smooth functions on $\Omega $ and the derivatives of all order of
these functions. Let $y\in {\cal A}_{\Omega }$ be a smooth function. Then
for $\lambda \in \Omega ,$ there exists a smooth function $z\in {\cal A}%
_{\Omega }$ such that $y=y\left( \lambda \right) e+\left( x-\lambda e\right)
z.$ Thus for each $\mu \in \Omega ,$ we have $\frac{dy}{dx}\left( \mu
\right) =z\left( \mu \right) e+\left( \mu -\lambda \right) \frac{dz}{dx}%
\left( \mu \right) .$ Since $z$ is smooth, there exists $\overline{z}\in
{\cal A}_{\Omega }$ such that $z=z\left( \lambda \right) e+\left( x-\lambda
e\right) \overline{z}.$ Therefore, $\frac{dy}{dx}\left( \mu \right) =z\left(
\lambda \right) +\left( \mu -\lambda \right) [\overline{z}\left( \mu \right)
+\frac{dz}{dx}\left( \mu \right) ].$ Hence $\frac{dy}{dx}=\frac{dy}{dx}%
\left( \lambda \right) +\left( x-\lambda e\right) \left( \overline{z}+\frac{%
dz}{dx}\right) .$ Now by induction on $n$ one can see that for each $n\geq 1$%
, there exists $y_{n}\in {\cal A}_{\Omega }$ such that $\frac{d^{n}y}{dx^{n}}%
=\frac{d^{n}y}{dx^{n}}\left( \lambda \right) e+\left( x-\lambda
e\right) y_{n}.$ By Remark 7.7, $\frac{d^{n}y}{dx^{n}}$ is a
smooth function. The rest of the proposition can be proved in the
same way.$\blacksquare $

{\bf Proposition\ 7.13. \ }Let $R$ be a semi-integral domain and $%
\mathbb{Q}
\subset R$ , or$R$ be an integral domain with characteristic zero. Let $%
\Omega \subset R$ be an absorbing subset. Assume that for some $\lambda \in
\Omega $ the ideal $\Lambda =I_{\lambda }$ of $C^{\infty }(\Omega )$ has the
property $\cap _{n=1}^{\infty }\Lambda ^{n}=\{0\}.$ Then the function $f\in
C^{\infty }(\Omega )$ is constant if and only if $\frac{df}{dx}=0.$

{\bf Proof. \ }If $f$ is constant then clearly $\frac{df}{dx}=0.$ Now assume
that $\frac{df}{dx}=0$. Then for each $n\in
\mathbb{N}
$ there exists $g_{n}\in C^{\infty }(\Omega )$ such that

\[
f-f(\lambda )e=(x-\lambda e)^{n}g_{n}.
\]%
So $f-f(\lambda )e\in \cap _{n=1}^{\infty }\Lambda ^{n}=0.$ Therefore $%
f=f(\lambda )e.\blacksquare $

A semi-integral domain is called {\it analytic }(resp. of {\it %
polynomial type}) if its smooth structure is analytic (resp. of
polynomial type). It is called {\it fine}, if there exists a
smooth function $\varphi :R-\{0\}\rightarrow R$ such that
$\varphi $ does not admit any smooth extension to $R$. A
semi-integral domain $R$ is called {\it wild }if
there exists a non-constant function $\psi :R\rightarrow R$ such that $\frac{%
d\psi }{dx}=0$. It is called {\it tame} if there exists a non-constant
smooth function $y:R\rightarrow R$ satisfying the following conditions:

There exists $\lambda \in R$ such that $\frac{d^{n}y}{dx^{n}}\left( \lambda
\right) =0,$\ for all\ $n=1,2,3,....$

{\bf Proposition\ 7.14. }\ Let $R$ be as in Proposition 7.13.
Then{\bf \ }$R$ is analytic if and only if it is not tame.

{\bf Proof. }\ Let{\bf \ }$R$ be analytic. Assume that $\varphi
:R\rightarrow R$ is smooth. Then, $\varphi \notin \cap
_{n=1}^{\infty }I_{0}{}^{n},$ where $I_{0}={\cal C}^{\infty
}(R)\cdot x.$ Therefore, there exists $n\in
\mathbb{N}
$ such that $\varphi \notin I_{0}{}^{n}$. In other words there exists no
function $\psi :R\rightarrow R$ such that $\varphi $ can be written as $\
x^{n}\cdot \psi $. Thus there exists $k\leq n$ such that $\frac{d^{k}\varphi
}{dx^{n}}\left( 0\right) \neq 0.$

The proof of the sufficiency is clear.$\blacksquare$

By the above proposition each smooth function on a analytic
semi-integral domain $R$ satisfying one of the conditions of
Proposition 7.13 is uniquely determined by its derivatives at an
element of $R$. Now assume that $y:R\rightarrow R$ is analytic.
Without any ambiguity we can write $y$ in the following form.

\[
\left( \ast \right) \ \ \ \ \ \ \ \ \ \ \ \ \ \ \ \ \ \ \ \ \ \ y=\Sigma
_{n=0}^{\infty }a_{n}\left( x-\lambda e\right) ^{n},\ \ \lambda \in R,
\]%
where \ $n!a_{n}=\frac{d^{n}y}{dx^{n}}\left( \lambda \right) $. The relation
$\left( \ast \right) $ is called the {\it series representation of} $y$ at $%
\lambda $. Moreover in this case we have:

{\bf Proposition\ 7.15. }\ The set of all the series
representations of smooth functions on $R$ at $\lambda \in R$, is
a commutative $R-$algebra under component-wise additions and
Cauchy products. Moreover, this algebra is closed under term-wise
differentiation and substitution.

{\bf Proposition}\ {\bf 7.16.} \ Let $X\ $be a non-coarse
connected topological space and let $R$ be the ring of all
continuous real functions on $X $. Then

1) For each $f\in C^{\infty }(%
\mathbb{R}
)$ the mapping $\overline{f}:R\rightarrow R$ given by $\overline{f}(\alpha
)=f\circ \alpha $ is smooth.

2) The mapping $f\longrightarrow \overline{f}$ $\ $\ from \ $C^{\infty }(%
\mathbb{R}
)$ into $C^{\infty }(R)$ given by $\overline{f}(\alpha )=f\circ \alpha $ is
a monomorphism of algebras.

3)\ Assume that $X=%
\mathbb{R}
.$ Then,$C^{\infty }(R)$ is wild.

{\bf Proof. \ }As we have seen earlier $R$ is a proper semi-integral domain.
So it has a unique smooth structure. Now for each smooth function $f\in
C^{\infty }(%
\mathbb{R}
)$ , and all $\alpha \in R$ define $\overline{f}(\alpha )=f\circ \alpha \in
R.$ By Lemma 5.5 there exists $g\in C^{\infty }(%
\mathbb{R}
)$ such that

\[
f(\beta (x))=f(\alpha (x))+(\beta -\alpha )(x)g(\beta (x)).
\]

Or

\[
\overline{f}(\beta )=\overline{f}(\alpha )+(\beta -\alpha )\overline{g}%
(\beta )
\]

By Remark 7.7, $\overline{f}$ is a smooth function on $R.$ The
rest of the proposition is immediate.$\blacksquare $

\section{ \ {\protect\Large Smooth Structure on }$%
\mathbb{Z}
$}

\bigskip\ \ \ \ {\bf Proposition\ 8.1. } Let{\bf \ }$\Omega \subset
\mathbb{Z}
$\ be an absorbing subset. Then the unique smooth structure on
$\Omega $ is analytic.

{\bf Proof. \ }Let $\omega $ be an element of $\Omega $ and let $f\in
I_{\omega }^{{\it n}}.$ Assume that for some $\lambda $ in $\Omega $ such
that $\lambda -\omega \neq \pm 1,$we have $f(\lambda )=\mu \neq 0.$ Then for
each $n\in
\mathbb{N}
$ , \ $(\lambda -\omega )^{n}\mid \mu .$ Since $%
\mathbb{Z}
$ is a unique factorization domain this is impossible.$\blacksquare $

{\bf Proposition\ 8.2. \ }Let $\ z:%
\mathbb{Z}
-\{{0}\}\rightarrow
\mathbb{Z}
$ be defined as follows

\[
z\left( t\right) =\left( 1-t^{2}\right) +2\left( 1-t^{4}\right) ^{2}\left(
2^{4}-t^{4}\right) ^{2}+2\left( 1+2\times 2^{8}\right) \left(
1-t^{18}\right) ^{3}\left( 2^{18}-t^{18}\right) ^{3}\left(
3^{18}-t^{18}\right) ^{3}+
\]

\[
...2\left( 1+2\times 2^{8}\right) \left( 1+2\times 2^{2\times 3^{3}}\times
3^{2\times 3^{3}}\right) \times ...\times \left( 1+2\Pi _{k=2}^{n}k^{2\times
n^{n}}\right) \times
\]

\[
\Pi _{k=1}^{n+1}(k^{2\times \left( n+1\right) ^{n}}-t^{2\times \left(
n+1\right) ^{n}})^{n+1}+....
\]

Then $z$ is a smooth function and does not admit any smooth extension to $%
\mathbb{Z}
$.

{\bf Proof. }$\ $Clearly $z$ is smooth. Assume that $\overline{z}:%
\mathbb{Z}
\rightarrow
\mathbb{Z}
$ is a smooth extension of $z$. Then, there exists a smooth function $y:%
\mathbb{Z}
\rightarrow
\mathbb{Z}
$ such that $\overline{z}=\overline{z}\left( 0\right) e+xy.$ Since $%
\overline{z}\left( 2\right) =\overline{z}\left( 0\right) +2y\left( 2\right)
=-3,$ $\overline{z}\left( 0\right) \neq 0$. Let

$I_{n+1}=\left( 1-t^{2}\right) +2\left( 1-t^{4}\right) ^{2}\left(
2^{4}-t^{4}\right) ^{2}+...+$

$2\left( 1+2\times2^{2\times2^{2}}\right) \left(
1+2\times2^{2\times3^{3}}\times3^{2\times3^{3}}\right) \times...\times\left(
1+\Pi_{k=2}^{n}k^{2\times n^{2}}\right) \times$

$\Pi_{k=1}^{n+1}\left( k^{2\left( n+1\right) ^{n}}-t^{2\left( n+1\right)
^{n}}\right) ^{n+1}$

and $\omega _{n}=\left( 1+2\times 2^{8}\right) \left( 1+2\times 2^{54}\times
3^{54}\right) \times ...\times \left( 1+2\Pi _{k=2}^{n}k^{2\times
n^{n}}\right) $. Clearly, $\omega _{n}$ divides $\overline{z}-I_{n+1}.$
Moreover, $\overline{z}-I_{n+1}$ is a smooth extension of $z-I_{n+1}$. There
exists $y_{1}\in {\cal A}_{%
\mathbb{Z}
}$ such that $\overline{z}-I_{n+1}=\left( \overline{z}-I_{n+1}\right) \left(
0\right) e+xy_{1}$. Therefore, $\left( \overline{z}-I_{n+1}\right) \left(
\omega _{n}\right) =\left( \overline{z}-I_{n+1}\right) \left( 0\right)
+\omega _{n}y_{1}\left( \omega _{n}\right) $. But $\left( \overline{z}%
-I_{n+1}\right) \left( \omega _{n}\right) =\left( z-I_{n+1}\right) \left(
\omega _{n}\right) $ is divisible by $\omega _{n}.$ Thus, $\left( \overline{z%
}-I_{n+1}\right) \left( 0\right) $ is divisible by $\omega _{n}.$
Furthermore, $I_{n+1}\left( 0\right) $ is divisible by $\omega _{n}$.
Therefore, for all $n\geq 3,$ $\overline{z}\left( 0\right) $ is divisible by
$\omega _{n}.$ This is a contradiction.$\blacksquare $

\bigskip As we have proved $%
\mathbb{Z}
$ is an analytic integral domain. The above proposition shows that $%
\mathbb{Z}
$ is also a fine integral domain. Some other properties of the smooth
functions on $%
\mathbb{Z}
$ is contained in $\left[ 2\right] $ and $[3].$ More about the
subject of this paper will be given later.

\begin{quote}
{\large Some Open Problems}
\end{quote}

1) Is there any $R-$algebra admitting non-consistent separated smooth
structures?

Let $({\cal A},\Sigma )$ be a smooth pair. A sub-module $M\subset {\cal A}$
is called a {\it basic sub-module} if for each $\Lambda \in \Sigma $, we
have $\pi _{\Lambda }(X)\in \Lambda ^{\circ }.$ It is clear that any two
basic sub-modules of $({\cal A},\Sigma )$ are isomorphic.

2) Is there any separated smooth $R-$algebra without any basic sub-module?

3) {\it Conjecture} : Every separated smooth $%
\mathbb{Z}
-$algebra is analytic.

4) Is every smooth function on $%
\mathbb{Z}
$ in the form given in Theorem 7.8 with $a_{i}\in{\mathbb{Z}}[x].$

5) Is there any semi-integral domain of polynomial type?

6) Characterize analytic,wild and fine semi-integral domains.

7) Let $\left( a_{n}\right) _{n\in
\mathbb{N}
}$ be a sequence in $%
\mathbb{Z}
$. Under what conditions does there exists a smooth function $y:%
\mathbb{Z}
\rightarrow
\mathbb{Z}
$ having $\Sigma _{n=0}^{\infty }a_{n}x^{n}$ as its series representation
around zero?

8) Characterize semi- integral domains $R$ having the property
that ${\cal C}^{\infty }(R)$ are semi- integral domains.

9) Is $(C^{\omega }\left(
\mathbb{R}
\right) ,[%
\mathbb{R}
^{{\it n}}])$ analytically and separatedly maximal?

10) Are any two analytically and separatedly maximal analytic
subalgebras of the $%
\mathbb{R}
-$algebra $%
\mathbb{R}
^{%
\mathbb{R}
}$ isomorphic?

11) Let $R$ be a semi-integral domain and let $\Omega \subset
\Omega
^{\prime }\subset R$ be absorbing subsets. Under what condition a function $%
h\in C^{\infty }(\Omega )$ can be extended to $\Omega ^{\prime }.$

\bigskip\ \ \ \ {\bf \ \ Acknowledgements}{\large :}\ Some parts of the work
has been done during the periods the author was at ICTP as an
associate member. He would like to thank them for their
hospitality. Lemma 7.1. is partly due to Prof. Rahim
Zaare-Nahandi. The author would like to thank him for this and
for reading the final version of this paper. He also thanks the
University of Tehran.

{\large References}

[1] Morris W. Hirsch , Differential Topology: Springer Verlag 1976

[2] A. Shafiei Deh Abad , On the theory of smooth structures :
ICTP, Trieste, Preprint IC/92/8

[3] A. Shafei Deh Abad, An introduction to the theory of
differentiable structures on infinite integral domains which are
not fields, J. Sc. I. R. Iran, Vol.1. No.3 (1990)

[4] A. Shafei Deh Abad, On the theory of smooth structures II ,
IC/92/268.

[5] A. Shafie Deh Abad and F. Kamali Khamseh, On the theory of
smooth structures III, IC/93/308. \ \ \ \ \ .

 [6]\ A. Shafei Deh Abad, On the theory of smooth structures IV,
IC/94/265.

\end{document}